

\baselineskip=14pt
\parskip=10pt

\font\eightrm=cmr8 
\font\eighttt=cmtt8
\magnification=\magstephalf

\def\1{{\overline{1}}}
\def\2{{\overline{2}}}
\parindent=0pt
\overfullrule=0in
\def\Tilde{\char126\relax}
\def\frac#1#2{{#1 \over #2}}
\bf
\centerline
{
Balls in Boxes:
}
\centerline
{
Variations on  a Theme of Warren Ewens and Herbert Wilf 
}
\rm
\bigskip
\centerline{ {\it
Shalosh B. EKHAD and Doron ZEILBERGER}\footnote{$^1$}
{\eightrm  \raggedright
Department of Mathematics, Rutgers University (New Brunswick),
Hill Center-Busch Campus, 110 Frelinghuysen Rd., Piscataway,
NJ 08854-8019, USA.
{\eighttt zeilberg  at math dot rutgers dot edu} ,
\hfill \break
{\eighttt http://www.math.rutgers.edu/\~{}zeilberg/} .
Written: June 13, 2011.
Accompanied by Maple package 
\hfill \break
{\eighttt http://www.math.rutgers.edu/\~{}zeilberg/tokhniot/BallsInBoxes} .
Sample input and output can be gotten from:
{\eighttt http://www.math.rutgers.edu/\~{}zeilberg/mamarim/mamarimhtml/bib.html} .
\hfill \break
Supported in part by the National Science Foundation of the United States of America.
}
}

\qquad \qquad\qquad \qquad
{\it To Herbert Saul Wilf (b. June 13, 1931), on his 80-th birthday}

{\bf Preface}

There are $r$ boys and $n$ girls. Each boy must pick {\it one} girl to invite to be his date in the prom.
Although each girl expects to get $R:=r/n$ invitations,
most likely, many of them would receive  less, and many of them would receive more. 
Suppose that Nilini, the most ``popular'' girl,
got as many as $m+1$ prom-invitations, is she indeed so popular, or did
she just ``luck-out''?

Each one of $r$ students has to choose from $n$ different parallel Calculus sections,
taught by different professors.
Although each  professor expects to get $R:=r/n$ students signing-up,
most likely, many of them would receive  less, and many of them would receive more. 
Suppose that Prof. Niles, the most ``popular'' professor
got as many as $m+1$ students, is Prof. Niles justified in assuming that
she is more popular than her peers, or did  she just ``luck-out''?

It is Saturday night, and there are $r$ people who have to decide where to dine,
and they have $n$ restaurants to choose from.
Although each restaurant expects to get $R:=r/n$ diners,
most likely, many of them would receive  less, and many of them would receive more. 
Suppose that the Nevada Diner, the most ``popular'' restaurant,
got as many as $m+1$ diners, can they congratulate themselves for the quality of
their food, or ambiance, or location, or can they only congratulate themselves for
being lucky?

Each one of $r$ cases of acute lymphocitic leukemia
has to choose one of $n$ towns (artificially made all with  equal-populations)
where to happen.
Although each town expects to get $R:=r/n$ cases,
most likely, many of them would receive  less, and many of them would receive more. 
Suppose that the Illinois town Niles had
$m+1$ cases of that disease, do its people have to
be concerned about their environment, or
is it only Lady Luck's fault?

Of course all these questions have the same answer, and typically one
talks about $r$ balls being placed, uniformly at random,
in $n$ boxes, where the largest number of balls that landed
at the same box was $m+1$.
Yet another way: A monkey is typing an $r$-letter
word using a keyboard of an alphabet with $n$ letters,
and the most frequent letter showed-up $m+1$ times.
Does the typing monkey have a particular fondness for that letter,
or is he a truly uniformly-at-random monkey who does not play favorites with the letters?

{\bf Asking the Right Question}

As Herb Wilf  pointed out so eloquently in his wonderful talk at
the conference W80 (celebrating his 80th birthday)
(based, in part, on [EW]), using the depressing disease formulation,
the right questions are {\bf not}:

``What is the probability that   Nilini would get so many ($m+1$ of them) prom-invitations?''

``What is the probability that Prof. Niles would get so many ($m+1$ of them) students?''

``What is the probability that the Nevada Diner would get so many ($m+1$ of them) diners?''

``What is the probability that Niles, IL would get so many ($m+1$ of them) cases of acute lymphocitic leukemia?''

Even though this is the wrong question (whose answer 
would make Nilini, Prof. Niles and the Nevada Diner's successes
go to their heads, and would make the real-estate prices in Niles, IL, plummet), because it is so tiny,
and seemingly extremely unlikely to be ``due to chance'', let's answer this question anyway.

The {\it a priori} probability of Nilini getting $m+1$ or more prom-invitations, using the
{\it Poisson Approximation} is:
$$
e^{-R} (\sum_{i=m+1}^{\infty} \frac{R^i}{i!})=
e^{-R} (e^R-\sum_{i=0}^{m} \frac{R^i}{i!})=
1-e^{-R}\sum_{i=0}^{m} \frac{R^i}{i!} \quad,
$$
indeed very small if $m$ is considerably larger than $R$.

But {\it a priori} we don't know who would be the ``lucky champion'' (or the unlucky town),
the {\bf right} question to ask is:

{\bf The Right Question:} Given $r$, $n$,  and $m$,
compute (if possible exactly, but at least approximately):

$P(r,n,m):=$: the probability that {\it every } box got $\leq m$ balls.

{\bf Getting the Right Answer to the Right Question, as Fast as Possible}

In [EW], Ewens and Wilf present a beautiful, {\it fast} ($O(mn)$), algorithm for
computing the {\it exact} value of $P(r,n,m)$, that employs a method that is described in the 
Nijenhuis-Wilf classic [NW] (but that has been around for a long time,
and rediscovered several times, e.g. by one of us ([Z1]), and before
that by J.C.P. Miller, and according to Don  Knuth the method goes back to Euler. At any 
rate, [EW]  does not claim novelty for the method, only for {\it applying} it to the present problem).

The {\it specific} real-life examples given in [EW] were:

1. (Niles, IL): $r=14400,n=9000$, (so $R=8/5$), $m=7$. Using their method, they got (in less than one second!) the value
$$
P(14400,9000,7)=0.0953959131671303999971555481626\dots \quad ,
$$
meaning that the probability that {\it every} town in the US, of the size of Niles, IL, would get
no more than $7$ cases is less than ten percent. So with probability $0.904604086832869600002844451837$,
{\it some} town (of the same size, assuming, artificially that the US has been divided into
towns of that size) somewhere, in the US, would get {\it at least} eight cases.  There is (most probably)
nothing wrong with their water, or their air-quality,  the only one that they may blame is Lady Luck!

For comparison, the {\it a priori} probability that Niles,IL would get $8$ or more cases is
roughly:
$$
1-e^{-1.6} \sum_{i=0}^{7} \frac{1.6^i}{i!}= 0.00026044\dots \quad ,
$$
a real reason for (unjustified!) concern.

2. (Churchill County, NV): $r=8000, n=12000$, (so $R=2/3$), $m=11$. Using their method, they got (in less than one second!) the value
$$
P(8000,12000,11)=0.999999895529647647310726013392\dots \quad ,
$$
so it is extremely likely that {\it every} district got at most $11$ cases, and the
probability that {\it some} district got $12$ or more cases is indeed small,
namely 
$$
1-P(8000,12000,11)=0.104470 \cdot 10^{-6} \quad ,
$$
so these people should  indeed panic.

For comparison, the {\it a priori} probability that Churchill County, NV, would get $12$ or more cases is
roughly:
$$
1-e^{-2/3} \sum_{i=0}^{11} \frac{(2/3)^i}{i!}=.870586315 \cdot 10^{-11} \quad,
$$
in that case people would have been right to be concerned, but for the wrong reason!

{\bf The Maple package } {\tt BallsInBoxes}

This article is accompanied by the Maple package {\tt BallsInBoxes}
available from: \hfill\break
{\tt http://www.math.rutgers.edu/\~{}zeilberg/tokhniot/BallsInBoxes} .

Lots of sample input and output files can be gotten from: \hfill\break
{\tt http://www.math.rutgers.edu/\Tilde zeilberg/mamarim/mamarimhtml/bib.html} .

\vfill\eject

{\bf How to Compute} $P(r,n,m)$ {\bf Exactly?}

Easy! As Ewens and Wilf point out in [EW], and Herb Wilf mentioned in his talk,
there is an obvious, explicit,  ``answer''
$$
P(r,n,m)=
\frac{1}{n^{r}}\sum \frac{r!}{r_1! r_2! \dots r_n!} \quad,
$$
where the sum ranges over the set of $n$-tuples of integers
$$
A(r,n,m):=\{(r_1, r_2, \dots, r_n) \,\, \vert \,\,
0 \leq r_1, \dots, r_n \leq m \quad , \quad r_1+r_2+ \dots +r_n=r \} \quad .
$$
So ``all'' we need, in order to get the {\it exact} answer,
is to construct the set $A(r,n,m)$ and add-up all the multinomial coefficients.

Of course, there is a better way. As is well-known (see [EW]), and easy to see,
writing
$$
P(r,n,m)=
\frac{r!}{n^r}\sum_{(r_1, \dots, r_n) \in A(r,n,m)} \frac{1}{r_1! r_2! \dots r_n!} \quad,
$$
the $\sum$ is the coefficient of $x^{r}$ in the expansion of 
$$
\left ( \sum_{i=0}^{m} \frac{x^i}{i!} \right )^n \quad ,
$$
so all we need is to go to Maple, and type (once $r,n$, and $m$ have been assigned numerical values)

{\tt r!/n**r*coeff(add(x**i/i!,i=0..m)**n,x,r);} \quad .

This works well for small $n$ and $r$, but, please, {\bf don't even try} to apply it to the first case of [EW],
($r=14400, n=9000, m=7$), Maple would crash!

Ewens and Wilf's brilliant idea was to use the Euler-Miller-(Nijenhuis-Wilf)-Zeilberger-\dots ``quick'' method for expanding a power of
a polynomial, and get an {\it answer} in less than a second!

[We implemented this method in Procedure {\tt Prnm(r,n,m)} of {\tt BallsInBoxes}].

While their method indeed takes less than a second (in Maple) for $r=14400, n=9000$ (and $7 \leq m \leq 12$), it takes quite a bit longer for
$r=144000, n=90000$, and I am willing to bet that for $r=10^8,n=10^8$ it would be hopeless to get
an {\it exact answer}, {\it even} with this fast algorithm.

But why this obsession with {\it exact} answers? Hello, this is {\it applied} mathematics, and
the epidemiological data is, of course, {\it approximate} to begin with, and we make lots of
unrealistic assumptions (e.g. that the US is divided into 9000 towns, each exactly the size of Niles, IL.) .
All we need to know is, ``are that many diseases likely to be due to pure chance, or is it a cause for concern'',
{\it Yes or No?}, {\it Ja oder Nein?}, {\it Oui ou Non?}, {\it Ken o Lo?}.

\vfill\eject

{\bf Enumeration Digression}

It would be nice to get a more compact (than the huge multisum above) (symbolic) ``answer'', or  ``formula'', in terms of
the {\it symbols} $r,n$ and $m$. This seems to be hopeless. But fixing, 
positive integers $a,b$ and $m$, one can ask for a ``formula'' (or whatever), in $n$,
for the quantity $P(an,bn,m)$ that can be written as $B(a,b,m;n)/(an)^{bn}$ where
$$
B(a,b,m;n):=
(an)!\sum_{(r_1, \dots, r_n) \in A(an,bn;m)} \frac{1}{r_1! r_2! \dots r_n!} \quad ,
$$
the cardinality of the {\it natural} combinatorial set consisting of
placing $an$ balls in $bn$ boxes in such a way that no box receives more than $m$ balls.
Equivalently, all {\it words} in a $bn$-letter alphabet, of length $an$, where no letter
occurs more than $m$ times. For example, when $a=b=m=1$, we have the deep theorem:
$$
B(1,1,1;n)=n! \quad .
$$
Equivalently, $e(n)=B(1,1,1;n)$ is a solution of the {\it linear recurrence equation with polynomial coefficients}
$$
e(n+1)-(n+1)e(n)=0 \quad, \quad (n \geq 0)  \quad,
$$
subject to the {\it initial condition}  $e(0)=1$.

It turns out that, thanks to the not-as-famous-as-it-should-be {\it Almkvist-Zeilberger} algorithm
[AZ] (an important component of the deservedly famous {\it Wilf-Zeilberger Algorithmic Proof Theory}),
one can find similar recurrences (albeit of higher order, so it is no longer ``closed-form'', in $n$)
for the sequences $B(a,b,m;n)$ for any {\it fixed} triple of positive integers, $a,b,m$.

(See Procedures {\tt Recabm} and {\tt RacabmV} in the Maple package {\tt BallsInBoxes}).

Indeed, since $B(a,b,m;n)$ is $(an)!$ times the coefficient of $x^{an}$ in
$$
\left ( \sum_{i=0}^{m} \frac{x^i}{i!} \right )^{bn} \quad ,
$$
it can be expressed, (thanks to {\it Cauchy}), as
$$
\frac{(an)!}{2\pi i}\oint_{|z|=1} \frac{\left ( \sum_{i=0}^{m} \frac{z^i}{i!} \right )^{bn}}{z^{an+1}} \, dz,
\eqno(Cauchy)
$$
and this is game for the Almkvist-Zeilberger algorithm, that has been incorporated into
\hfill\break
{\tt BallsInBoxes}. See the web-book

{\tt http://www.math.rutgers.edu/\Tilde zeilberg/tokhniot/oBallsInBoxes2} 

for these recurrences for $1 \leq a, b \leq 3$ and $ 1 \leq m \leq 6$.

{\bf Asymptotics}

Once the first-named author of the present article computed a recurrence,
it can go on, thanks to the {\it Birkhoff-Trzcinski method} ([WimZ][Z2]),
to get very good asymptotics! So now we can get a very precise 
asymptotic formula (in $n$) (to any desired order!)
for $P(an,bn,m)$, that turns out to be very good for large, and even not-so-large $n$, and for {\it any}
desired $a,b,m$. 
Procedure {\tt Asyabm} in our Maple package {\tt BallsInBoxes}
finds such asymptotic formulas. See

{\tt http://www.math.rutgers.edu/\Tilde zeilberg/tokhniot/oBallsInBoxes1}

for asymptotic formulas, derived by combining Almkvist-Zeilberger with {\tt AsyRec} (also included in
{\tt BallsInBoxes} in order to make the latter self-contained.)

This works for {\it every} $m$, and {\it every} $a$ and $b$, in principle!
In practice, as $m$ gets larger than $10$, the recurrences become very high order, and
take a very long time to derive.

But as long as $m \leq 8$ and even (in fact, especially) when $n$ is very large,
this method is much faster than the method of [EW]  ($O(mn)$ with large $n$ is not that small!).
Granted, it does not give you an {\it exact} answer, but neither do they (in spite of their claim,
see below!) .

But let's be pragmatic and forget about our purity and obsession with ``exact'' answers.
Since we know from ``general nonsense'' that the asymptotics of 
the desired probability 
$$
C(a,b,m;n):=P(an,bn,m) \quad (=B(a,b,m;n)/(an)^{bn})
$$
behaves asymptotically as
$$
C(a,b,m; n) \asymp \mu^n(c_0+ O(1/n) )\quad,
$$
for {\it some}  numbers $\mu$ and $c_0$,
all we have to do is crank out (e.g.) the $200$-th and $201$-th term 
and estimate $\mu$ to be $C(a,b,m;201)/C(a,b,m;200)$, and then
estimate $c_0$ to be $C(a,b,m;200)/\mu^{200}$. Using Least Squares one can do even better,
and also estimate higher order asymptotics (but we don't bother, enough is enough!).

Procedure {\tt AsyabmEmpir} in our Maple package {\tt BallsInBoxes}
uses this method, and gets very good results!

For example, for the Niles, IL, example, in order to get estimates for $P(14400,9000,m)$, typing

{\tt evalf(subs(n=1800,AsyabmEmpir(8,5,m,200,n)));}

for $m=7,8,9,10,11,12$ yields (almost  instantaneously)

$m=7$: $ 0.09540287131\dots$ (the exact value being: $0.095395913167\dots$) ,

$m=8$: $0.664971462304\dots$ (the exact value being: $0.66495441\dots$) ,

$m=9$: $0.9378712268719\dots$ (the exact value being: $0.93786433\dots$) ,

$m=10$: $0.990845139\dots$ (the exact value being: $0.9908433\dots$) ,

$m=11$: $0.998789295\dots$ (the exact value being: $0.99878892861\dots$) .

The advantage of the present approach is that we can handle very large $n$, for example, with the
same effort we can compute

{\tt evalf(subs(n=180000,AsyabmEmpir(8,5,m,200,n)))}

getting, for example,  that $P(1440000,900000,11)$ is  very close to $0.88554890636027$. The method used in [EW] 
(i.e. typing {\tt Prnm(1440000,900000,11);} in {\tt BallsInBoxes}) would take forever!

{\bf Caveat Emptor}

There is another problem with the $O(mn)$  method described in [EW]. Sure enough, it works well for the
examples given there, namely $P(14400,9000,m)$ for $6 \leq m \leq 12$ and
$P(8000,12000,m)$ for $4 \leq m \leq 8$.

This is corroborated by our implementation of that method, 
(Procedure  {\tt Prnm(r,n,m)} in \hfill\break {\tt BallsInBoxes}).

Typing  (once {\tt BallsInBoxes} has been read onto a Maple session):

{\tt t0:=time(): Prnm(14400,9000,9) , time()-t0; }

returns

$0.937864339305858219725360911354, 0.884$

that tells you the desired value (we set {\tt Digits} to be $30$), and that it took $0.884$ seconds to compute that value.

But now try:

{\tt t0:=time(): Prnm(1000,100,15), time()-t0;}

and get in $0.108$ seconds (real fast!)
 
$-0.728465229161818857989128673465 \cdot 10^{50}$ \quad .

Something is rotten in Denmark! We learned in kindergarten that a {\it probability} has to be between $0$ and $1$, 
so a negative probability, especially one with 50 decimal digits, is a bit fishy.
Of course, the problem is that [EW]'s ``exact'' result is not really {\it exact}, as it uses floating-point arithmetic.

Big deal, since we work in Maple, let's increase the system variable {\tt Digits} (the number of digits used in
floating-point calculations), and type the following line:

{\tt evalf(Prnm(1000,100,15),80);}

getting $5.71860506564981...$,
a little bit better! (the probability is now less than six, and at least it is positive!), but still nonsense.

{\tt Digits:=83} still gives you nonsense, and it only starts to ``behave'' at {\tt Digits:=90}.

Now let's multiply the inputs, $r$ and $n$ by $10$, and take $m=22$ and try to evaluate
$P(10000,1000,22)$. Even {\tt Digits:=250} still gives nonsense! Only {Digits:=310} gives you something
reasonable and (hopefully) correct.

The way to overcome this problem is to keep upping {\tt Digits} until you get close
answers with both {\tt Digits} and, say, {\tt Digits+100}. This is implemented in Procedure
\hfill\break
{\tt PrnmReliable(r,n,m,k)} in {\tt BallsInBoxes},
if one desires an accuracy of $k$ decimal digits. This is {\it reliable}
indeed, but {\bf not} exact, and {\it not} rigorous, since it uses numerical heuristics.
The exact answer is a {\it rational number}, that is implemented in Procedure {\tt PrnmExact(r,n,m)} of
{\tt BallsInBoxes}.

{\bf The Cost of Exactness}

If you type

{\tt  t0:=time():PrnmExact(14400,9000,7): time()-t0;}

you would get in $42$ seconds (no longer that fast!) a {\it rational number} whose numerator and denominator
are {\it exact} integers with $54207$ digits.

See 

{\tt http://www.math.rutgers.edu/\Tilde zeilberg/tokhniot/oBallsInBoxes7a} 

for the outputs (and timings) of {\tt PrnmExact(14400,9000,m);} for $m$ between $6$ and $12$
and see 

{\tt http://www.math.rutgers.edu/\Tilde zeilberg/tokhniot/oBallsInBoxes7b}

for the outputs (and timings) of {\tt PrnmExact(8000,12000,m);} for $m$ between $4$ and $8$.
No longer fast at all! ($2535$ and $248$ seconds respectively).

{\bf Let's Keep It Simple: An Ode to the Poisson Approximation}

At the end of [EW], the authors state:

{\it `` A Poisson Approximation is also possible but it may be  inaccurate, particularly around the tails
of the distribution. Our exact method is fast and does not suffer from any of those problems.''}

Being curious, we tried it out, to see if it is indeed so bad. Surprise, it is terrific!
But let's first review the  Poisson approximation as we understand it.

The probability of any particular box (of the $n$ boxes) getting $\leq m$ ball is, roughly, using the Poisson approximation
($R:=r/n$):
$$
e^{-R}  \sum_{i=0}^{m} \frac{R^i}{i!} \quad .
$$
Of course the $n$ events are {\bf not} independent, but let's pretend that they are. The probability that
{\it every} box got $\leq m$ balls is approximated by
$$
Q(r,n,m):=
\left ( e^{-R}  \sum_{i=0}^{m} \frac{R^i}{i!} \right )^n \quad .
$$

[$Q(r,n,m)$ is implemented by procedure {\tt PrnmPA(r,n,m)} in
{\tt BallsInBoxes}. It is as fast as lightning!]

Ewens and Wilf are very right when they claim that $P(r,n,m)$ and $Q(r,n,m)$ are very far apart
around the ``tail'' of the distribution, but who cares about the tail? Definitely not a scientist
and even not an applied mathematician. It turns out, empirically (and we did extensive 
numerical testing,
see Procedure {\tt HowGoodPA1(R0,N0,Incr,M0,m,eps)} in {\tt BallsInBoxes}), that whenever
$P(r,n,m)$ is not extremely small, it is very well approximated by $Q(r,n,m)$, and using
the latter (it is so much faster!) gives very good approximations, and enables one to construct
the ``center'' of the probability distribution (i.e. ignoring the tails) very accurately.
See 

{\tt http://www.math.rutgers.edu/\Tilde zeilberg/tokhniot/oBallsInBoxes4} \quad ,

and

{\tt http://www.math.rutgers.edu/\Tilde zeilberg/tokhniot/oBallsInBoxes5} \quad ,

for comparisons (and timings!, the Poisson Approximation wins!) .

In particular, the estimates for the {\it expectation}, {\it standard deviation}, and even the higher moments
match extremely well!

Another (empirical!) proof of the fitness of the Poisson Approximation can be seen in:

{\tt http://www.math.rutgers.edu/\Tilde zeilberg/tokhniot/oBallsInBoxes1}

where the (rigorous!) asymptotic formulas derived, via {\tt AsyRec}, from the recurrences obtained via the
Almkvist-Zeilberger algorithm are very close to those predicted by the Poisson Approximation (except for very small $m$,
corresponding to the ``tail'').

{\bf The Full Probability Distribution of the Random Variable ``Maximum Number of Balls in the Same Box''}

It would be useful, for  given positive integers $a$ and $b$, to know how the probability distribution
``maximum number of balls in the same box when throwing $an$ balls into $bn$ boxes'' behaves. 
One can ``empirically'' construct (without arbitrarily improbable tail) the distribution 
of the random variable ``maximum number of balls in the same box'' when $an$ balls are
uniformly-at-random placed in $bn$ boxes
(Let's call it $X_n(a,b)$, and $X_n$ for short)
using
$$
Pr(X_n=m)=P(an,bn,m)-P(an,bn,m-1) \quad .
$$
First, and foremost, what is the expectation, $\mu_n$?
Second what is the standard deviation, $\sigma_n$?, skewness?, kurtosis?, and it would ve even nice to know
higher  $\alpha$-coefficients (alias moments of $Z_n:=(X_n-\mu_n)/\sigma_n$), 
as asymptotic formulas in $n$.

For the expectation, $\mu_n$, Procedure {\tt AveFormula(a,b,n,d,L,k)}
uses the more accurate ``empirical approach'' and  Maple's built-in 
Least-Squares command,
to obtain the following empirical (symbolic!) estimates for the expectation.

$a=1,b=1$: {\tt evalf(AveFormula(1,1,n,1,300,1000,10),10);} yields that

$\mu_n$ is roughly  $2.293850526 + (0.4735983525 ) \cdot \log n $

$a=2,b=1$: {\tt evalf(AveFormula(2,1,n,1,300,1000,10),10);} yields that

$\mu_n$ is roughly  $3.963420618 + (0.5834252496) \cdot \log n$

$a=1,b=2$: {\tt evalf(AveFormula(1,2,n,1,300,1000,10),10);} yields that

$\mu_n$ is roughly  $1.640094145 + (0.3873602232) \cdot \log n$.

Note that for $a=1,b=1$, the approximation to $\mu_n$ can be written 
\hfill\break
$2.293850526 + (1.090500507) \cdot \log_{10} n$,
so a ``rule-of-thumb'' estimate for the expectation when $n$ balls are thrown into $n$ boxes is
a bit more than $2$ plus the number of (decimal) digits.

Procedure {\tt NuskhaPA1(R,n,K,d)} uses the Poisson Approximation to guess polynomials in $\log n$ of degree $d$ fitting
the average, standard deviation, and higher moments, as 
asymptotic expressions in $n$, for
$nR$ balls thrown into $n$ boxes, where $R$ is now any (numeric) {\it rational} number.
Even $d=1$ seems to give a fairly good fit, so they all seem to be (roughly) linear in $\log n$.

{\bf Procedure} {\tt SmallestmPA}

Procedure
{\tt SmallestmPA(r,n,conf)} gives you the smallest $m$ for which, with confidence {\tt conf}, you can deduce
that the high value of $m$ is {\bf not} due to chance (using the Poisson Approximation). For example

{\tt SmallestmPA(14400,9000,.99);}

yields $10$, meaning that if a town the size of Niles, IL got $10$ or more cases, then with probability $>0.99$ 
it is not just bad luck. If you want to be $\%99.99$-sure of being a victim of the environment rather
than of Lady Luck, type:

{\tt SmallestmPA(14400,9000,.9999);}

and get $13$, meaning that if you had $13$ cases, then with probability larger than $0.9999$ it is not due to chance.

{\bf The Minimum Number of Balls that Landed in the Same Box, Procedure } {\tt LargestmPA}

An equally interesting, and harder to compute, random variable is the {\it minimum number of balls that landed in the same box},
but the Poisson Approximation handles it equally well. Analogous to {\tt SmallestmPA}, we have, in {\tt BallsInBoxes},
Procedure {\tt LargestmPA(r,n,conf)} that
tells you the largest $m$ for which you can't blame luck for getting $m$ or less balls.

For example, if there are $10000$ students that have to decide between $100$ different calculus sections,

{\tt LargestmPA(10000,100,.99); }

that happens to be $66$, tells you that any section that only has $66$  students or less, with probability $>0.99$,
it is because that professor (or time slot, e.g. if it is an 8:00am class) is not popular, and you can't blame bad luck.

{\tt  LargestmPA(10000,100,.9999);}

that outputs $57$, tells you that anyone who only had $ \leq 57$ students enrolled is unpopular with probability $>\%99.99$,
and can't blame bad luck.

On the other end, going back to the original problem, \hfill\break
{\tt  SmallestmPA(10000,100,.99);}  yields $139$, telling you  that any section for which $139$ or more students signed up
is {\it probably} (with prob. $>0.99$) due to the popularity of that section, while
{\tt  SmallestmPA(10000,100,.9999);}  yields $151$.

{\bf Final Comments}

1. One can possibly (using the {\it saddle-point method})
get asymptotic formulas from the contour integral $(Cauchy)$, but this is not {\it our} cup-of-tea, so
we leave it to other people.

2. Another ``back-of-the-envelope'' ``Poisson Approximation'' is to argue that since the
probability of any individual box getting strictly more than $m$ balls is roughly (recall that $R=r/n$)
$$
e^{-R}  \sum_{i=m+1}^{\infty} \frac{R^i}{i!} =
e^{-R} (e^R- \sum_{i=0}^{m} \frac{R^i}{i!})=
1- e^{-R}\sum_{i=0}^{m} \frac{R^i}{i!} ,
$$
by the {\it linearity of expectation}, the expected number of {\it lucky} (or {\it unlucky} if the balls are diseases) boxes
exceeding $m$ balls is roughly
$$
n \left ( 1- e^{-R}\sum_{i=0}^{m} \frac{R^i}{i!} \right ) \quad .
$$
In the case of Niles, IL, the expected number of towns that would get $8$ or more cases is:
$$
9000 \left ( 1- e^{-1.6}\sum_{i=0}^{7} \frac{(1.6)^i}{i!} \right )= 2.343961376410372 \quad,
$$
so it is not at all surprising that at least one town got as many as $8$ cases. On the other hand,
in the other example $r=8000, n=12000,m=12$, the expected number of unfortunate counties is:
$$
12000 \left ( 1- e^{-(2/3)}\sum_{i=0}^{12} \frac{(2/3)^i}{i!} \right )= 0.533706802 \cdot 10^{-8} \quad,
$$
so it is indeed a reason for concern.

{\bf Conclusion}

We completely agree with Ewens and Wilf that {\it Monte Carlo} takes way too long, and is not that accurate, 
and that {\it their} method is far superior to it. But we strongly disagree with their
dismissal of the  Poisson Approximation.
In fact, we  used their ingenious method to conduct extensive empirical (numerical) testing that established that
the  Poisson Approximation, that they dismissed as 
``inaccurate'', is, as a matter of fact, sufficiently accurate,
far more reliable,
in addition to being yet-much-faster!
It is much safer to use the Poisson Approximation than to use their ``exact'' method
(in floating-point arithmetic), and when one uses {\it truly} exact calculations,
in rational arithmetic, their ``fast'' method becomes {\it anything but}.

Even when the floating-point problem is addressed by using multiple precision ({\tt PrnmReliable} discussed above),
their fast algorithm becomes slow for very large $r$ and $n$, while the Poisson Approximation is
almost  instantaneous even for very large $r$ and $n$, and {\it any} $m$.

So while we believe that the algorithm in [EW] is not as {\it useful} as the Poisson Approximation,
it sure was {\it meta-useful}, 
since it enabled us to conduct extensive numerical testing that showed,
{\it once and for all}, that it is far less useful then the latter.

Additional evidence comes from our own symbolic approach
(fully rigorous for $m \leq 9$ and semi-rigorous for higher values of $m$),
that establishes the 
adequacy of the Poisson Approximation for {\it symbolic} $n$.

Finally, as we have already pointed out,
since the data that one gets in applications is always {\it approximate} to begin with, insisting on an ``exact'' answer,
even when it is easy to compute, is  unnecessary.

{\bf Coda: But We, Enumerators, Do Care About Exact Results!}

Our point, in this article, was that for {\it applications} to
statistics, the Poisson Approximation suffices. But {\it we}
are {\it not} statisticians. We are {\it enumerators}, and
we do like exact results! The approach of [EW] enables us to know,
for example,
in less than one second the {\bf exact} number of ways that 
$1001$ balls can be placed in $1001$ boxes such that no box
received more than $7$ balls. Just type (in {\tt BallsInBoxes})

{\tt (1001**1001)*PrnmExact(1001,1001,7); }

and get a beautiful {\bf exact} integer with $3004$ digits!

Typing

{\tt (1001**1001)*PrnmPA(1001,1001,7); }

will give you something fairly close (the ratio being $0.9997852\dots$)

but  for a {\bf pure} enumerator, this is very unsatisfactory.
So long live exact answers!, but {\it not} in statistics.

{\bf References}

[AZ] Gert Almkvist and Doron Zeilberger, {\it The Method of Differentiating Under The Integral Sign},
J. Symbolic Computation {\bf 10} (1990), 571-591. 
[Available on-line from: 
\hfill\break
{\tt http://www.math.rutgers.edu/\Tilde zeilberg/mamarim/mamarimPDF/duis.pdf}]

[EW] Warren J. Ewens and Herbert S. Wilf,
{\it Computing the distribution of the maximum in balls-and-boxes problems with application to clusters of disease cases},
Proc. National Academy of Science (USA) {\bf104(7)} (July 3, 2007), 11189-1191 .
[Available on-line from: 
\hfill\break
{\tt http://www.pnas.org/content/104/27/11189.full.pdf}]

[WimZ] Jet Wimp and Doron Zeilberger,
{\it Resurrecting the asymptotics of linear recurrences},
J. Math. Anal. Appl. {\bf 111} (1985), 162-177.
[Available on-line from: 
\hfill\break
{\tt http://www.math.rutgers.edu/\Tilde zeilberg/mamarimY/WimpZeilberger1985.pdf}]

[Z1] Doron Zeilberger, {\it The J.C.P. Miller Recurrence for exponentiating a polynomial, and its q-Analog},
J. Difference Eqs. and Appl. {\bf 1} (1995), 57-60.
[Available on-line from:
\hfill\break
{\tt http://www.math.rutgers.edu/\Tilde zeilberg/mamarim/mamarimhtml/power.html}]

[Z2] Doron Zeilberger
{\it AsyRec: A Maple package for Computing the Asymptotics of Solutions of Linear Recurrence Equations with Polynomial Coefficients},
The Personal Journal of Shalosh B. Ekhad and Doron Zeilberger {\tt http://www.math.rutgers.edu/\Tilde zeilberg/pj.html},
April 6, 2008.
[Article and package available on-line from:
\hfill\break
{\tt  http://www.math.rutgers.edu/\Tilde zeilberg/mamarim/mamarimhtml/asy.html}]

\end